\documentclass[a4paper,12pt]{amsart}
\usepackage{amssymb}
\usepackage{ifthen}
\usepackage{graphicx}
\usepackage{float}
\usepackage{caption}
\usepackage{subcaption}
\usepackage{cite}
\usepackage{amsfonts}
\usepackage{amscd}
\usepackage{amsxtra}
\usepackage{color}

\usepackage[latin1]{inputenc}

\newtheorem{thm}{Theorem}
\newtheorem{cor}{Corollary}
\newtheorem{lem}{Lemma}
\newtheorem{rem}{Remark}
\newtheorem{prop}{Proposition}

\newtheorem*{conj}{Conjecture}
\newtheorem{prob}{Problem}

\theoremstyle{definition}
\newtheorem{defn}{Definition}[section]
\newtheorem{example}{Example}

\newenvironment{pf}[1][]{%
	\vskip 1mm
	\noindent
	\ifthenelse{\equal{#1}{}}%
	{{\slshape Proof. }}%
	{{\slshape #1.} }%
}%
{\qed\bigskip}

\newcounter{alphabet}



\def\be{\begin{equation}}
	\def\ee{\end{equation}}

\newcommand{\bee}{\begin{enumerate}}
	\newcommand{\eee}{\end{enumerate}}

\newcommand{\blem}{\begin{lem}}
	\newcommand{\elem}{\end{lem}}
\newcommand{\bthm}{\begin{thm}}
	\newcommand{\ethm}{\end{thm}}
\newcommand{\bcor}{\begin{cor}}
	\newcommand{\ecor}{\end{cor}}
\newcommand{\beg}{\begin{example}}
	\newcommand{\eeg}{\end{example}}
\newcommand{\begs}{\begin{examples}}
	\newcommand{\eegs}{\end{examples}}
\newcommand{\bdefe}{\begin{defn}}
	\newcommand{\edefe}{\end{defn}}
\newcommand{\bprob}{\begin{prob}}
	\newcommand{\eprob}{\end{prob}}
\newcommand{\bques}{\begin{ques}}
	\newcommand{\eques}{\end{ques}}
\newcommand{\bei}{\begin{itemize}}
	\newcommand{\eei}{\end{itemize}}
\newcommand{\bcon}{\begin{conj}}
	\newcommand{\econ}{\end{conj}}
\newcommand{\bcons}{\begin{conjs}}
	\newcommand{\econs}{\end{conjs}}
\newcommand{\bprop}{\begin{propo}}
	\newcommand{\eprop}{\end{propo}}
\newcommand{\br}{\begin{rem}}
	\newcommand{\er}{\end{rem}}
\newcommand{\brs}{\begin{rems}}
	\newcommand{\ers}{\end{rems}}
\newcommand{\bo}{\begin{obser}}
	\newcommand{\eo}{\end{obser}}
\newcommand{\bos}{\begin{obsers}}
	\newcommand{\eos}{\end{obsers}}
\newcommand{\bpf}{\begin{pf}}
	\newcommand{\epf}{\end{pf}}
\newcommand{\ba}{\begin{array}}
	\newcommand{\ea}{\end{array}}
\newcommand{\beq}{\begin{eqnarray}}
	\newcommand{\beqq}{\begin{eqnarray*}}
		\newcommand{\eeq}{\end{eqnarray}}
	\newcommand{\eeqq}{\end{eqnarray*}}

\newcounter{minutes}
\divide\time by 60
\newcounter{hours}
\multiply\time by 60 \addtocounter{minutes}{-\time}

\begin{document}
	\title[Tischler-Tyson inequality]{Tischler-Tyson type inequality for k-symmetric polynomials}

	\author{Ilgiz R. Kayumov}
	\address{St~Petersburg University, Universitetskaya Embankment 7-9, 199034, St~Petersburg, Russia \\
		Kazan Fedearl University, Kremlevskaya 18,
			420 008 Kazan, Russia}
	%\curraddr{}
	\email{i.kayumov@gmail.com}
	%\thanks{}
	
	%    author two information
	\author{Diana M. Khammatova}
	\address{St~Petersburg University, Universitetskaya Embankment 7-9, 199034, St~Petersburg, Russia; Keldysh Institute of Applied Mathematics RAS, Moscow, Russia}
	%\curraddr{}
	\thanks{Corresponding author: Diana Khammatova}
	\email{anakinlynx@gmail.com}

	

	\begin{abstract}
			In this paper, we prove a sharp inequality for the critical points and
			critical values of $k$-symmetric complex polynomials with a restriction on their critical points. The inequality was inspired by the papers of Tischler (1989) and Tyson (2005) and has previously been established only in a few special cases. As a corollary, our result implies Dual Smale's Mean Value conjecture with an improved constant for this class of polynomials.
	\end{abstract}

	\maketitle
	\pagestyle{myheadings}
	
	\markboth{I.R. Kayumov, D.M. Khammatova}{Tischler-Tyson inequality}
	
	\section{Introduction}
	
	The present paper continues a long-standing line of research on the critical points and critical values of complex polynomials, initiated by S. Smale in 1981. In his paper ~\cite{s}, he formulated the conjecture that is now known as Smale's Mean Value Conjecture.
	
	Let $\mathfrak{P}_n$ denote the class of polynomials $f$ satisfying the normalization conditions $f(0)=0$ and $f'(0)=1$. Smale's Mean Value Conjecture claims that for every $f\in\mathfrak{P}_n$, there is a critical point $\zeta$ such that
	$$
	\left|\frac{f(\zeta)}{\zeta}\right|\leq 1-\frac1n.
	$$
	
	For the polynomial $f(z)=z+cz^n$ the value of $\left|\frac{f(\zeta)}{\zeta}\right|$ is equal to $1-\frac1n$ at every critical point $\zeta$ of $f$, so the constant cannot be improved. 
	
	This conjecture attracted considerable attention and many attempts were made to prove or disprove it. Some estimates can be found in the papers \cite{bmn} and \cite{cfl}. Results for polynomials with real critical points were obtained in \cite{hk} and later generalized in \cite{d26}. Further results related to Smale's Mean Value Conjecture are discussed in Section 7.2 of the book \cite{rs}. 
	
	Recently, Smale's Mean Value Conjecture with the constant $1$ was disproved using AI tools. A geometric mechanism for constructing counterexamples is presented in a recent preprint by T. Adamczewski and M. J. Colbrook \cite{ac}. 
	
	One can also consider classes of functions other than polynomials. The papers \cite{ng17} and \cite{d24} contain results on finite Blaschke products. In the preprint \cite{jnoct}, finite Blaschke products are used to prove the existence of polynomials of arbitrarily large degree for which the corresponding minimum exceeds $1.1018$. In particular, this result also disproves the classical conjecture.

	In 2009, V.~N.~Dubinin and T.~Sugawa \cite{d}, and independently T.~W.~Ng \cite{ng}, introduced what is now known as Dual Smale's Mean Value Conjecture. It can be formulated as follows:

	\par\medskip
	
	\textbf{Dual Smale's Mean Value Conjecture.}
	\textit{
			For every $f\in\mathfrak{P}_n$, there exists a critical point $\zeta$ of $f$ satisfying
			$    \left|\frac{f(\zeta)}{\zeta}\right|\geq\frac1n$}
	
	Dubinin and Sugawa \cite{d} first proved that for all $f\in\mathfrak{P}_n$, $\max\limits_{f'(\zeta)=0}\left|\frac{f(\zeta)}{\zeta}\right|\geq\frac{1}{n4^n}$. Further developments in this direction are discussed in the review paper~\cite{d12}.

	An important result was obtained by Dubinin in 2019 \cite{dub19}. He showed that for every complex polynomial $f$ of degree $n \geq 2$ and every $z \in \mathbb{C}$, there exists a critical point $\zeta$ such that $f'(\zeta) = 0$ and
	$$
	\left|	\frac{f(z)-f(\zeta)}{z-\zeta}\right|\geq \frac1n\tan\frac{\pi}{4n}|f'(z)|.
	$$ 
	The author noted that the constant $\tan\frac{\pi}{4n}$ can be replaced by $\frac{1}{n}$, which yields an estimate of order $\frac{1}{n^2}$.
	
	In a recent preprint \cite{jnaug}, A. Jatar and T. W. Ng obtained an improved lower bound $\frac1{n^2}\left(n-\frac12\right)^{1/n}$ and established several other results concerning both the classical and dual Smale conjectures.
	
	While the general case remains open, the exact inequality is known in several special cases. A number of useful inequalities can be found in D. Tischler's 1989 paper \cite{t}. For instance, he proved the required inequality for conservative polynomials, i.e. those for which $f(\zeta)/\zeta$ has the same value at all critical points. 
	
	It is natural to consider the conjecture for small values of $n$. The cases $n\leq 4$ follow from Tischler's work. It turned out that as the value of $n$ increases, the amount of computation increases rapidly, making a direct verification of the conjecture impractical. 
	
	In the papers \cite{hkk19} and \cite{hkk25}, the value $f(\zeta)/\zeta$ was estimated at a critical point of minimal modulus. This approach made it possible to prove the conjecture for $n\leq 7$ but could not be extended to $n=8$, as  shown in \cite{hkk19}. This made it necessary to find another approach.

	A promising idea is suggested by Tischler's aforementioned paper \cite{t}. He proved the following inequality: 	
	\begin{equation*}
		\left|\frac{f(\zeta)}{\zeta}-\frac{1}{2}\right| \leq \frac{1}{2}-\frac{1}{n}, \quad 2 \leq n \leq 4.
	\end{equation*}
	
	This estimate seems very important, as it implies both the classical and dual conjectures. Although J. T. Tyson \cite{tyson} showed that it does not hold for $n\geq 5$, constructing a counterexample, Tischler's paper plays an important role in the current research. In the course of his proof, he obtained the following inequality:
	
	\begin{equation*}		\sum_{i=1}^3\left|\frac{\frac{f(\zeta_i)}{\zeta_i}-\frac12}{\zeta_i^2}\right|^2 \leq \frac1{16}	\sum_{i=1}^3\left|\frac1{\zeta_i^2}\right|^2.
	\end{equation*}
	
	It was later generalized by Tyson, who proved that
	$$ \sum_{j=1}^{n-1} \left|z_j^{2-n}\left(\frac{1}{2}-\frac{f(z_j)}{f'(0)z_j}\right)\right|^2 \leq
	K(n)^2 \left(\frac1{|z_1|^{2(n-2)}}+\ldots+\frac1{|z_{n-1}|^{2(n-2)}}\right)
	$$ for a constant $K(n)$ depending only on $n$.
	
	We call the inequalities of this form Tischler-Tyson type inequalities. Their main advantage is that they take all critical points into account simultaneously without having to select one.
	
	In \cite{hk26}, the following functional was introduced:	
		\begin{equation*}\label{ourhyp}
				F_n:=\left|\frac1{\zeta_1^{n-1}}P\left(\zeta_1\right)\right|^2+\ldots+\left| \frac1{\zeta_{n-1}^{n-1}}P\left(\zeta_{n-1}\right)\right|^2 - \frac1{n^2}\left(\frac1{|\zeta_1|^{2n-4}}+\ldots+\frac1{|\zeta_{n-1}|^{2n-4}}\right).
		\end{equation*}
	
	It was proved that this expression is non-negative in case when all critical points of a polynomial $f\in\mathfrak{P}_n$ lie on a single ray. This result immediately implies Dual Smale's Mean Value Conjecture for this class of polynomials. In \cite{hk26AMS}, the inequality $F_n\geq 0$ was proved for all polynomials in $\mathfrak{P}_n$ with $n\leq 5$. In the same paper, an example was constructed for $n=9$ showing that it does not hold in the general case. However, numerical evidence suggested that it may be true for polynomials with real critical points.
	  
	On the other hand, both this inequality and the Dual conjecture itself can be generalized. This idea was motivated by the preprint of Q. Tang \cite{tang}, who proved Dual Smale's Mean Value Conjecture for odd polynomials using the constant $\frac1{n^2}$ obtained by Dubinin in \cite{dub19}. He also noted that the proof can be extended to $k$-symmetric polynomials, i.e. polynomials of the form $f(z) = zQ(z^k)$, where $Q$ is a polynomial.
	
	The proof is short and clear, but the result is not sharp, as Tang himself points out. This suggests that a stronger constant should exist for $k$-symmetric polynomials.
	
	A candidate for the constant is the value $C(n,k) = \frac{1}k B(\frac1k,\frac{n-1}k+1)$. The reasoning behind this choice can be found in \cite{dk26}. If $k=1$, we get the class $\mathfrak{P}_n$ and $C(n,1)=\frac1n$.
	
	The main aim of the present paper is to prove the following theorem:
	\begin{thm}\label{mainth}
		Let $f(z)\in\mathfrak{P}_n$ be a $k$-symmetric polynomial with critical points $\zeta_1,\ldots,\zeta_{n-1}$, such that $(\zeta_j)^k\in\mathbb{R}$ for all $1\leq j\leq n-1$. Then 
		\begin{equation}\label{ksymm}
			\left|\frac{f\left(\zeta_1\right)}{\zeta_1^{n-k}} \right|^2+\ldots+\left| \frac{f\left(\zeta_{n-1}\right)}{\zeta_{n-1}^{n-k}}\right|^2 - C(n,k)^2\left(\frac1{|\zeta_1|^{2(n-k-1)}}+\ldots+\frac1{|\zeta_{n-1}|^{2(n-k-1)}}\right) \geq 0,
		\end{equation}
		where $C(n,k) = \frac{1}k B(\frac1k,\frac{n-1}k+1)$.	
	\end{thm}

The paper is organized as follows. Section 2 outlines the main idea of the proof and formulates a technical proposition that is crucial to completing it. Section 3 presents two lemmas necessary for the detailed proof given in Section 4. The final section summarizes the results obtained and explains their significance.

	\section{The main idea}
	
	A $k$-symmetric polynomial $f\in \mathfrak{P}_n$ of degree $n=k\cdot m+1$ can be represented as $f(z) = zQ(z^k)$, where $Q$ is a polynomial of degree $m$, $Q(0)=1$. Taking the derivative, we see that
	$$
	f'(z) = Q(z^k)+kz^kQ'(z^k),
	$$ 
	which together with the condition $(\zeta_j)^k\in\mathbb{R}$ gives the critical points $(r_j e^{2\pi i s})^{1/k}$, $0\leq s\leq k-1$, $1\leq j \leq m$. We denote $w_j=\frac1{r_j}$ and suppose that $|w_j|\leq 1$, which can always be achieved by a suitable normalization.
	
	Then, as shown in \cite{dk26}, we have the following representation:
	\begin{equation}\label{krepr}
		\frac {f(z)}z=\int\limits_0^1\left(1-t^k z^k w_1\right)\cdot\ldots\cdot \left(1-t^k z^k w_m\right)\,dt.
	\end{equation}
	
	As $f(z) = z Q(z^k)$, we see that
	$$
	\frac{f(\zeta_{(j-1)k+1})}{\zeta_{(j-1)k+1}} =\ldots = \frac{f(\zeta_{jk})}{\zeta_{jk}} = Q\left(\frac1{w_j}\right).
	$$ 
	
	Therefore, \eqref{ksymm} takes the form
	\begin{equation*}
		\left|w_1^{m-1} Q\left(\frac1{w_1}\right)\right|^2+\ldots+\left| w_m^{m-1}Q\left(\frac1{w_m}\right)\right|^2 - C(m,k)^2\left(|w_1|^{2(m-1)}+\ldots+|w_{m}|^{2(m-1)}\right) \geq 0,
	\end{equation*}
	where $C(m,k) = \frac{1}k B(\frac1k,m+1)$.
	
	If $m=1$, this expression is equal to zero and the proof is trivial. The case $m=2$ was proved in \cite{dk26}. In what follows, we take $m\geq 3$. 
	
	We can simplify the left-hand side of this inequality by subtracting the following non-negative quantity:
	$$
	w_1^{2(m-1)}\left|Q\left(\frac1{w_1}\right)-C(m,k)\right|^2+\ldots + w_m^{2(m-1)}\left|Q\left(\frac1{w_m}\right)-C(m,k)\right|^2
	$$	
	
	The remaining part is
	
	\begin{equation*}\label{comexp}
		2C(m,k)\sum\limits_{j=1}^m \left(Q\left(\frac1{w_j}\right)-C(m,k)\right) w_j^{2(m-1)},
	\end{equation*}		 
		
	Expanding this expression and dividing it by $2k C(m,k)$, we get
	
	\begin{multline}\label{posexp}
		G(w_1,\ldots,w_m):=\sum\limits_{j=1}^m \frac{(-1)^{j+1}}{(kj-(k-1))(kj+1)}\sum\limits_{i_s\in[1,m]} \left( \prod\limits_{s=1}^j w_{i_s}\sum_{s=1}^j w_{i_s}^{2m-2-j}\right) - \\ - \frac{C(m,k)}k \sum\limits_{j=1}^m w_j^{2m-2},	
	\end{multline}
	where the sum $\sum\limits_{i_s\in[1,m]}$ runs over all subsets of $j$ distinct indices from ${1,\ldots,m}$.

	To proceed, we need a technical but important inequality, which we prove below.
	
	\begin{prop}\label{mainprop}
	Let $M,N$ be integers with $M\geq 1$ and $N\geq M+\max(\theta-2,1)$, where $\theta>0$, and let $w_j\in[0,1]$ for $1\leq j\leq M$. Then

\begin{equation*}
		g := \sum\limits_{j=1}^{M} \frac{(-1)^{j+1}}{(j+\theta-1)(j+\theta)}\sum\limits_{i_s\in[1,M]}\left( \prod\limits_{s=1}^j w_{i_s}\sum_{s=1}^j w_{i_s}^{N-j}\right) - B\left(\theta,M+1\right)\sum\limits_{j=1}^{M} w_j^{N} \geq 0.
\end{equation*}
	
	\end{prop}
	
	First, let $w_j>0$ for $1\leq j\leq m$. Then Proposition~\ref{mainprop} with $M=m$, $N=2m-2$ and $\theta=\frac1k$ yields the non-negativity of \eqref{posexp}. If $k=1$, we obtain the desired inequality for polynomials from $\mathfrak{P}_n$ with positive critical points; this case was proved in \cite{hk26}. The case $k=2$ corresponds to odd polynomials with real critical points.
	
	To prove the general case, we apply the algorithm outlined in \cite{hk26AMS} for polynomials in $\mathfrak{P}_9$. Let $w_1,\ldots, w_{m-u}$ be positive and $w_{m-u+1},\ldots,w_m$ be negative. For $j=1,\ldots, u$ we set $v_j = -w_{j+m-u}$. Thus, we obtain two sets of positive variables: $w_1,\ldots, w_{m-u}$ and $v_1,\ldots, v_u$. We know that
	
	$$
	G(w_1,\ldots, w_{m-u}, v_1,\ldots, v_u) \geq 0.
	$$
	
	Let us consider the difference
	$$
	G(w_1,\ldots, w_{m-u},  w_{m-u+1},\ldots, w_{m}) - G(w_1,\ldots, w_{m-u}, v_1,\ldots, v_u).
	$$
	
	Using \eqref{posexp}, we obtain the following expression:
	\begin{multline*}
		\sum\limits_{c=1}^u\sum\limits_{j=1}^{m-u} \frac{(-1)^{j+1}(1-(-1)^c)}{(k(j+c)-(k-1))(k(j+c)+1)} \sum\limits_{i_s\in[1,m-u]} \left( \prod\limits_{s=1}^j w_{i_s}\sum_{s=1}^j w_{i_s}^{N_1(j,c)}\right) \sum\limits_{i_s\in[1,u]}  \prod\limits_{s=1}^c v_{i_s} + \\ + \sum\limits_{c=1}^u\sum\limits_{j=1}^{m-u} \frac{(-1)^{c+1}(1-(-1)^j)}{(k(j+c)-(k-1))(k(j+c)+1)} \sum\limits_{i_s\in[1,u]} \left( \prod\limits_{s=1}^c v_{i_s}\sum_{s=1}^c v_{i_s}^{N_1(j,c)}\right) \sum\limits_{i_s\in[1,m-u]}  \prod\limits_{s=1}^j w_{i_s},
	\end{multline*}
	where $N_1(j,c) = 2m-2-j-c$.
	
	Since all the variables are positive, it suffices to prove that 
	$$
	\sum\limits_{j=1}^{M_1} \frac{(-1)^{j+1}}{(k(j+c)-(k-1))(k(j+c)+1)} \sum\limits_{i_s\in[1,M_1]} \left( \prod\limits_{s=1}^j w_{i_s}\sum_{s=1}^j w_{i_s}^{N_1(j,c)}\right) \geq 0
	$$
	for all $k\geq 1$, $1\leq c \leq m-1$, $1\leq M_1\leq m-c$. We take
	$M_1 = m-c$, as the results for smaller $M_1$ follow by setting some of the variables equal to zero. This gives
	$$
		\frac1{k^2}\sum\limits_{j=1}^{m-c} \frac{(-1)^{j+1}}{(j+c+\frac1k-1)(j+c+\frac1k)} \sum\limits_{i_s\in[1,m-c]} \left( \prod\limits_{s=1}^j w_{i_s}\sum_{s=1}^j w_{i_s}^{2m-2-j-c}\right).
	$$
	
	Applying Proposition~\ref{mainprop} with $M=m-c$, $N=2m-2-c$ and $\theta = c+\frac1k$, we conclude that this expression is bounded below by
	$$
	\frac1{k^2}B\left(c+\frac1k,m-c+1\right)\sum\limits_{j=1}^{m-c}w_j^{2m-2-c},
	$$
	which is positive for all positive $w_j$.
	
	Then,
	$$
	G(w_1,\ldots, w_{m-u},  w_{m-u+1},\ldots, w_{m}) > G(w_1,\ldots, w_{m-u}, v_1,\ldots, v_u)\geq 0.
	$$
	
	Therefore, if we prove Proposition~\ref{mainprop}, Theorem~\ref{ksymm} will hold.

%	\section{Small values of $m$}
%	{\color{blue} Add a section for $m=1,2$?}
	
	\section{Necessary lemmas}
	\begin{lem}\label{lem1}
		For integers $p\geq 1$, $t\geq 0$, $N\geq t+p+1$ and for $x\in [0,1]$
		$$
		\frac{(N-p-1)!}{N!}\sum\limits_{j=0}^t (-1)^j\binom{t}{j}\frac{(N-j)!}{(N-j-p-1)!}x^j\geq(1-x)^{t}.
		$$
	\end{lem}
	\begin{proof}
		Let $F(x) =	\frac{(N-p-1)!}{N!}\sum\limits_{j=0}^t (-1)^j\binom{t}{j}\frac{(N-j)!}{(N-j-p-1)!}x^j$. 
		
		Setting $y=1-x$, we write
		$$
		F(1-y) = \sum\limits_{j=0}^t c_j y^j,
		$$
		where
		$$
		c_j =\frac{(N-p-1)!}{N!}\binom{t}{j}\frac{(N-j)!}{(N-j-p-1)!} {}_2F_1(1-N+p+j,j-t,j-N;1).
		$$
		
		The prefactor of the hypergeometric function is evidently positive. The function can be represented this way:
		$$
		{}_2F_1(1-N+p+j,j-t,j-N;1) = \frac{(p+1)\cdot p\cdot\ldots\cdot(p+2-t+j)}{(N-j)\cdot(N-j-1)\cdot\ldots\cdots(N-t+1)}
		$$
		This expression is positive for $j>t-p-2$ and equals zero for $j\leq t-p-2$.
		
		Thus, $c_j\geq 0$ for all $j=0,\ldots,t$. Since $F(0)=1$, we have
		$$
		1=F(0)=\sum\limits_{j=0}^t c_j.
		$$
		
		Then
		$$
		F(x) = F(1-y) = \sum\limits_{j=0}^tc_j y^j\geq y^t \sum\limits_{j=0}^t c_j = y^t = (1-x)^t.
		$$
		Lemma~\ref{lem1} is proved.
	\end{proof}
	
	\begin{lem}\label{lem2}
		For integers $p\geq 1$, $2\leq b\leq p$, $t\geq 0$ and $N\geq p+\max(\theta-1,t+1)$ and for $\theta>0$ 
		\begin{multline*}
			D:=\int\limits_0^1 \Biggl[ \frac{p!}{(p-b)!}(-1)^b x^{b+\theta-1}(1-x)^{p-b+t+1}+ px^{\theta-1}(1-x)^pH_{N,b}(x)- \\ - px^{\theta}(1-x)^pH_{N-1,b}(x) - (p+1-b)x^{\theta-1}(1-x)^{p+1}H_{N,b}(x)	
			\Biggr]\, dx - \\ - \int\limits_0^1(b-1)b!\binom{N}{N-b}x^{\theta-1}(1-x)^{t+p+1}\, dx > 0,
		\end{multline*}
		where 
		$$
		H_{N,b}(x) := \sum\limits_{j=0}^t\frac{(-1)^j(N-j)!}{(N-j-b)!}\binom{t}{j} x^j.
		$$
		
	\end{lem}
	
	\begin{proof}
		We set $y=1-x$ and use Vandermonde's identity to rewrite $	H_{N,b}(x)$.
		
		\begin{equation}\label{hnest}
			H_{N,b}(x)  = b!\sum\limits_{j=0}^{\min(b,t)}\binom{N-t}{b-j}\binom{t}{j}y^{t-j}\geq b!y^{t}\sum\limits_{j=0}^{\min(b,t)}\binom{N-t}{b-j}\binom{t}{j} = b!y^{t}\binom{N}{b}.
		\end{equation}
		
		Using the recurrence relation $H_{N,b}(x) = H_{N-1,b}(x)+bH_{N-1,b-1}(x)$, we obtain
		
		\begin{equation*}
			H_{N,b}(x) - xH_{N-1,b}(x) = yH_{N-1,b}(x)+bH_{N-1,b-1}(x).
		\end{equation*}

	 	This allows us to rewrite part of the integrand as follows:
		\begin{multline*}
			x^{\theta-1} y^p(H_{N,b}(x)(p-(p+1-b)(1-x))-H_{N-1,b}(x)px) = \\ =x^{\theta-1} y^p((b-1)H_{N,b}(x)+x((p+1-b)(H_{N-1,b}(x)+b H_{N-1,b-1}(x))-pH_{N-1,b}(x))) = \\ = 	
			x^{\theta-1} y^p\Bigl((b-1)yH_{N-1,b}(x)+(b(b-1)+xb(p+1-b))H_{N-1,b-1}(x)\Bigr).
		\end{multline*}
		
		Since all the coefficients are positive, we can apply the estimate \eqref{hnest}. Estimating the first term from below by replacing $(-1)^b$ by $-1$, we obtain		
		\begin{multline*}
			D\geq-\frac{p!}{(p-b)!}B\left(b+\theta,p-b+t+2\right) +  (b-1)b!\binom{N-1}{b}B\left(\theta,p+t+2\right) + \\ +
			(b-1)b!\binom{N-1}{b-1}B\left(\theta,p+t+1\right) + 
			 b!(p+1-b)\binom{N-1}{b-1}B\left(\theta+1,p+t+1\right) - \\ -
			(b-1)b!\binom{N}{N-b}B\left(\theta,t+p+2\right).
		\end{multline*}
		
		Using the identity
		\begin{equation*}
			B\left(\theta,p+t+1\right) - B\left(\theta,p+t+2\right) = \int\limits_0^1x^{\theta} (1-x)^{p+t}\, dx = B\left(\theta+1,p+t+1\right),
		\end{equation*}
		
		we finally obtain	
		\begin{multline*}
			D\geq -\frac{p!}{(p-b)!}B\left(b+\theta,p-b+t+2\right) -
			(b-1)b!\binom{N}{N-b}B\left(\theta,p+t+2\right) +\\
			+(b-1)b!\left(\binom{N-1}{b}+\binom{N-1}{b-1}\right)B\left(\theta,p+t+2\right) +\\ +
			(b-1)b!\binom{N-1}{b-1}B\left(\theta+1,p+t+1\right) +
			b!(p+1-b)\binom{N-1}{b-1}B\left(\theta+1,p+t+1\right) = \\ = 
		b!\left(p\binom{N-1}{b-1}B\left(\theta+1,p+t+1\right) -\binom{p}{b}B\left(b+\theta,p-b+t+2\right)\right).
		\end{multline*}
		
	Then, $D\geq b!(H_1 - H_2)$, where 		
	\begin{multline*}
			\frac{H_1}{H_2} = \frac{p\binom{N-1}{b-1}\Gamma\left(p+t+\theta+2\right)\Gamma\left(\theta+1\right)\Gamma\left(p+t+1\right)}{\binom{p}{b}\Gamma\left(b+\theta\right)\Gamma\left(p-b+t+2\right)\Gamma\left(p+t+\theta+2\right)} = \\ = \frac{b\binom{N-1}{b-1}(p+t)\cdot\ldots\cdot (p+t-b+2)}{\binom{p-1}{b-1}\left(b+\theta-1\right)\cdot\ldots\cdots(1+\theta)}.
	\end{multline*}
	
	If $\theta<1$, we immediately have $H_1> H_2$, as $p> b+\theta-1$ and $N > p$. Otherwise, we use  $N\geq \theta+b-1$ and $t\geq 0$. Then 
	\begin{equation*}
		\frac{H_1}{H_2} \geq\frac{b \Gamma(b+\theta-1)\Gamma(p-b+1)\Gamma(\theta+1)\Gamma(p+1)}{\Gamma(\theta)\Gamma(p)\Gamma(b+\theta)\Gamma(p-b+2)} = \frac{b\theta p}{(p-b+1)(b+\theta-1)}.
	\end{equation*}
	Since $b\theta p- (p-b+1)(b+\theta-1) = (b-1)(b-1+p(\theta-1)+\theta)> 0$ for $\theta\geq 1$, then $H_1- H_2> 0$. Therefore, $D> 0$. 	This completes the proof.

	\end{proof}
	
	\section{Proof of Proposition~\ref{mainprop}}
	
	The present proof generalizes the algorithm introduced in \cite{hk26}. We introduce the following notation.
	$$
	\alpha(j):=\frac{(-1)^{j+1}}{(j+\theta-1)(j+\theta)}, \qquad \lambda:=  B\left(\theta,M+1\right).
	$$	
%	Here $\alpha(j)$ are the coefficients from the expression $g$ and $\lambda$ plays the same role as the constant in the Theorem.
	In what follows, we will work with two expressions	
	$$
	S_1(q,j_1,j_2) :=  \sum\limits_{j=j_1}^{j_2} \alpha(j)\sum\limits_{i_s\in[1,j_2-q]}w_{i_1}\ldots w_{i_{j-q}}(w_{i_1}^{N-j}+\ldots + w_{i_{j-q}}^{N-j}),
	$$
	and
	$$
	S_2(q,j_1,j_2,\phi(j)) := \sum\limits_{j=j_1}^{j_2}\phi(j) \sum\limits_{i_s\in[1,j_2-q]}w_{i_1}\ldots w_{i_{j-q}}.
	$$
	
	We start with the expression
	$$
	g_0^0:=g = S_1(0,1,M)-\lambda \sum\limits_{j=1}^{M} w_j^N.
	$$
	
	Let $d>0$. We define the operator $T_d(p)$ acting on polynomials $p(x_1,\ldots, x_m)$ by
	 $$
	 T_d(p):= d\cdot p - \left(x_1\frac{\partial p}{\partial x_1} +\ldots + x_m\frac{\partial p}{\partial x_m}\right).
	 $$
	
	This operator has several useful properties. First, it eliminates all homogeneous terms of degree $d$. In particular, if the degree of $p$ is equal to $d$, then the degree of $T_d(p)$ is at most $d-1$. For a homogeneous polynomial $p$ of degree $d$, we have $T_d(p)=0$. Moreover, if $T_d(p)$ is non-negative on a compact set, then $p$ is non-negative at every interior extremum. Consequently, to establish the non-negativity of $p$ throughout the set, it suffices to check its values on the boundary.

	The polynomial $g_0^0$ is homogeneous, its degree is equal to $N$. Applying the operator $T_N$, we obtain $0$ and check the values on the boundary of the set $[0,1]\times \cdots \times [0,1]$, which mean setting one of the variables equal to $0$ or $1$. To get the complete proof, we repeat the following two operations:	
	
	\begin{enumerate}
		\item Set $p$ of the variables equal to $1$. We denote 
		$$
		g_0^p:= g_0^0\Big|_{w_{M-p+1}=\ldots = w_{M}=1}.
		$$
		A direct calculation shows that $g_0^M=0$, so it suffices to consider $p\leq M-1$.
		
		\item Apply the operator repeatedly to reduce the degree of $g_0^p$ until the resulting polynomial becomes easier to investigate. We denote
		$$
		g_b^p = T_{N-b+1}(g_{b-1}^p).
		$$
		
		Since $N-b+1$ is the degree of $g_{b-1}^p$, applying $T_{N-b+1}$ reduces the degree by at least one. Consequently, $g_N^p$ is a constant, and we need only consider $b\leq N$.

	\end{enumerate}	
	
	The algorithm can be described as follows. For each fixed $p$ we assume that the non-negativity of $g_{b}^{p+1}$ has already been established. This is always possible for $p=M-1$, since $g_b^M=0$. Then we take $g_0^p$ and successively apply the corresponding operators until the non-negativity of $g_{b_0}^p$ can be established for some $b_0$. Having done this, we return to the polynomials $g_b^p$ with $b<b_0$ and check their values on the boundary. 
	
	Thus, the sufficient conditions for $g$ to be non-negative are
	
	\begin{enumerate}
		\item $g_{b_0}^p\geq0$ for some $b_0\in \mathbb{N}$;
		
		\item $h_b^p:=g_{b}^p\Big|_{w_{M-p}=1}\geq 0$ for all $1\leq b<b_0$;
		
		\item $u_b^p:=g_{b}^p\Big|_{w_{M-p}=0}\geq 0$ for all $0\leq b<b_0$.
	\end{enumerate}
	
	We note that $h_0^p = g_0^{p+1}$, so they do not have to be checked separately.
	
	We now work with $S_1$ and $S_2$ and examine how the operator $T_d$ acts on these two expressions. Calculations show that
		
	\begin{equation}\label{tds1}
		T_d S_1(q,j_1,j_2) = (d+q-N) S_1(q,j_1,j_2),
	\end{equation}
	
	\begin{equation}\label{tds2}
		T_d S_2(q,j_1,j_2,\phi(j)) = S_2(q,j_1,j_2,\phi(j)(d-j+q)) .
	\end{equation}

	Now we calculate $g_0^p$. We set $p$ variables in $g_0^0$ equal to $1$ and obtain
	\begin{multline*}
		g_0^p = \sum\limits_{s=0}^p\Biggl( \binom{p}{s} S_1(s,s+1,M-p+s) + s\binom{p}{s}S_2\left(s,s,M-p+s,\alpha(j)\right)\Biggr)   - \lambda p - \lambda \sum\limits_{j=1}^{M-p} w_j^N = \\ =
		\sum\limits_{s=0}^p \binom{p}{s} S_1(s,s+1,M-p+s)+S_2\left(0,0,M-p,\gamma(j,p)\right) - \lambda p - \lambda \sum\limits_{j=1}^{M-p} w_j^N,
	\end{multline*}
	where $\gamma(j,p):=(-1)^jp B\left(j+\theta,p+1\right)$.
	
	To get $g_b^p$, we apply the operators $T_{N},\ldots,T_{N-b+1}$. From \eqref{tds1} and \eqref{tds2},
	\begin{multline*}
		g_b^p = \sum_{s=b}^p \frac{s!}{(s-b)!}\binom{p}{s} S_1(s,s+1,M-p+s) + \\ + S_2\left(0,0,M-p,\gamma(j)\cdot\frac{(N-j)!}{(N-j-b)!}\right)-\lambda p \frac{N!}{(N-b)!}.
	\end{multline*}
	
	For $b=p+1$ we get
	$$
	g_{p+1}^p = S_2\left(0,0,M-p,\gamma(j)\cdot\frac{(N-j)!}{(N-j-p-1)!}\right)-\lambda p \frac{N!}{(N-p-1)!}.
	$$
	
	This expression is linear with respect to each variable. Therefore, its minimum is attained at a vertex, where each $w_j$ is either $0$ or $1$. Since the expression is symmetric, its value at a vertex depends only on the number of coordinates equal to $1$. Let $W_t$ denote the value of $g_{p+1}^p$ when $w_1=\ldots=w_t=1$ and all the other $w_j=0$, $0\leq t\leq M-p\leq N-p-1$. Then
	\begin{multline*}
		W_t = p\sum\limits_{j=0}^t\binom{t}{j}\frac{(-1)^j (N-j)!}{(N-j-p-1)!} B\left(j+\theta,p+1\right) - p\frac{N!}{(N-p-1)!} B\left(\theta,M+1\right) = \\ =
	p \int\limits_0^1 \frac{N! x^{\theta-1}(1-x)^p}{(N-p-1)!}  \left(\frac{(N-p-1)!}{N!}\sum\limits_{j=0}^t \binom{t}{j}\frac{(-1)^j(N-j)!}{(N-j-p-1)!}x^j-(1-x)^{M-p}\right)\,dx.
	\end{multline*}
		
	By Lemma~\ref{lem1} and the inequality $t\leq M-p$ we see that $W_t\geq 0$. It means that $g_{p+1}^p\geq 0$ for all admissible values of the variables. The first sufficient condition is satisfied with $b_0=p+1$.
	
	We next calculate $h_b^p$. A direct computation gives
	\begin{multline*}
		h_b^p = \frac{p+1-b}{p+1}\sum_{s=b}^{p+1}\frac{s!}{(s-b)!}\binom{p+1}{s} S_1(s,s+1,M-p+s-1)+ \\ + S_2\Bigl(0,0,M-p-1,\frac{p!}{(p-b+1)!}\gamma(j+b,p-b+1) + \frac{(N-j)!}{(N-j-b)!} \gamma(j,p) + \\ +
		\frac{(N-j-1)!}{(N-j-b-1)!} \gamma(j+1,p) \Bigr) - \lambda p\frac{N!}{(N-b)!}.
	\end{multline*}
	
	The term involving $S_1$ coincides with the corresponding term in $g_b^{p+1}$,
	which is non-negative by assumption. Thus, it suffices to prove that
	$$
	q_b^p := (p+1) h_b^p-(p+1-b)g_b^{p+1}\geq 0.
	$$
	
	Again, we obtain an expression that is linear in each variable: 
	\begin{multline*}
		q_b^p = S_2\Biggl(0,0,M-p-1,\frac{(p+1)!}{(p-b+1)!}\gamma(j+b,p-b+1)+ \\ + (p+1)\frac{(N-j)!}{(N-j-b)!} \gamma(j,p) +(p+1)\frac{(N-j-1)!}{(N-j-b-1)!} \gamma(j+1,p) - \\ - (p+1-b)\gamma(j,p+1)\frac{(N-j)!}{(N-j-b)!} \Biggr) - \lambda(p+1)(b-1)\frac{N!}{(N-b)!}.
	\end{multline*}
	
	For $b=1$ this expression equals zero. Let $b\geq 2$.
	
	As before, we set $t$ variables equal to $1$ and the remaining variables equal to $0$ for $0\leq t\leq M-p-1$. We want to prove that 
	$$
	D := (p+1)(J_0+J_1+J_2+J_3+J_4)\geq 0,
	$$
	where 
	$$
	J_0 =  -\int\limits_0^1(b-1)b!\binom{N}{N-b}x^{\theta-1}(1-x)^M\,dx,
	$$
	$$
		J_1 =  
		\int\limits_0^1\frac{p!}{(p-b)!}(-1)^b x^{b+\theta-1}(1-x)^{p-b+t+1} \, dx,
	\qquad
	J_2 = 
	\int\limits_0^1 p x^{\theta-1}(1-x)^p
		H_{N,b}(x) \, dx,
	$$
	$$
	J_3 = -\int\limits_0^1p x^{\theta}(1-x)^{p}H_{N-1,b}(x)\,dx, \qquad
	J_4 = -\int\limits_0^1(p+1-b) x^{\theta-1}(1-x)^{p+1}H_{N,b}(x)\,dx. 
	$$
	Here $H_{N,b}$ is defined as in Lemma~\ref{lem2}.

	Since $t+p+1\leq M<N$ and $p+\theta-1\leq M+\theta-2\leq N$, all the hypotheses of Lemma~\ref{lem2} are satisfied, and hence $D\geq0$. Consequently, $h_b^p\geq q_b^p\geq0$ for $b<b_0$. Thus, the second sufficient condition is satisfied.
	
	It remains to check that $u_b^p\geq 0$ for $1\leq p\leq M-1$ and $0\leq b\leq p$. The inequality is already established for $b=p+1$. We fix $b$ and $p$ and suppose that $u_{b+1}^p\geq 0$.  
	
	The polynomial $u_b^p$ has degree $N-b$. We apply the operator and see that
	$$
	T_{N-b}\left(g_b^p\Big|_{w_{M-p} = 0}\right)=(T_{N-b}g_b^p)\Big|_{w_{M-p}=0} = u_{b+1}^p\geq 0.
	$$ 
	
	The second sufficient condition gives $u_b^p\Big|_{w_{M-p-1}=1}=h_b^p\Big|_{w_{M-p-1}=0}\geq 0$. It remains to check $u_b^p\Big|_{w_{M-p-1}=0}$. Repeating the same
		steps, we come to the expression
	$$
	g_b^p\Big|_{w_1=\ldots = w_{M-p} = 0} = \frac{p N!}{(N-b)!}\left(B\left(\theta,p+1\right)-B\left(\theta,M+1\right)\right).
	$$
	
	Evidently, this expression is non-negative for $p\leq M$. The third sufficient condition holds, which finishes the proof of Proposition~\ref{mainprop}
	
	\section{Conclusion}
	
In this paper, we prove a generalization of the Tischler-Tyson type inequality for $k$-symmetric polynomials  proposed in \cite{hk26}. This inequality first appeared in \cite{dk26}, where it was proved  for $n=3k+1$. For $k=1$, the result was established in \cite{hk26} for polynomials in $\mathfrak{P}_n$ with positive critical points. For $k=2$, it was proved under the condition  $(\zeta_j)^k>0$, which corresponds to odd polynomials with real critical points.  Theorem~\ref{mainth} covers all these cases, allowing  $(\zeta_j)^k$ to take arbitrary real values, and establishes the inequality for all normalized $k$-symmetric polynomials whose critical points satisfy this condition.

Therefore, the present theorem  generalizes several previous results by establishing the desired estimate for $k$-symmetric polynomials of arbitrary degree. It also improves the constant in Dual Smale's Mean Value Conjecture for this class of polynomials. Although the Tischler-Tyson type inequality is known to fail for polynomials with arbitrary complex critical points, the improved constant may still be valid in the context of the Dual Smale conjecture for $k$-symmetric polynomials. This possibility opens a new direction for further research.

%	The present paper is devoted to a counterpart of the Tischler-Tyson type inequality for $k$-symmetric polynomials. In \cite{dk26}, this inequality was proved for $m=3$ and arbitrary $k$. Here, we fix $k=2$ and prove it for all $m$, under an additional technical restriction on the critical points.
%	
%	The obtained result can be regarded as a generalization of \cite{hk26} to odd polynomials with real critical points. At the same time, it provides a counterpart of Dual Smale's Mean Value Conjecture for this class of polynomials, with a better constant $C_m$.
%	
%	The proof strategy used in \cite{hk26} in the present paper was modified to work for the new conditions. It is possible that the same approach can be taken to prove a counterpart of Theorem~\ref{mainth} for the common case of $k$-symmetric polynomials. Eventually, there is hope to prove the inequality $G_{m,k}\geq 0$ for all real $\{\zeta_j^k\}$.
%	
%	It is also important to note that the constant $C(m,k)$ is a significant improvement compared  to $\frac1n$ in the dual conjecture. One can check that $\lim\limits_{k\to\infty} C(m,k) = 1$, hence for each fixed value of $m$ we can find $k$ for which this value will be close to $1$.

	\section*{Acknowledgments}
	The work of I. Kayumov was supported by Russian Science
	Foundation (grant 23-11-00066-P)

	\section*{Conflicts of interest}
	The authors have no conflicts of interest to declare that are relevant to the content of this article.
	
	\section*{Data availability statement}
	No datasets were generated or analysed during the current study. All logical and mathematical steps are self-contained within the manuscript.

\end{document}